\documentclass[draft]{amsart}
\usepackage{amssymb,amscd}
\usepackage[all]{xy}

\newcommand{\bH}{\mathbf H}
\newcommand{\bL}{\mathbf L}

\newcommand{\Z}{\mathbb Z}
\newcommand{\Q}{\mathbb Q}

\newcommand{\C}{\mathbb C}
\newcommand{\N}{\mathbb N}

\newcommand{\cA}{\mathcal A}

\newcommand{\cE}{\mathcal E}
\newcommand{\cF}{\mathcal F}
\newcommand{\cG}{\mathcal G}
\newcommand{\cH}{\mathcal H}
\newcommand{\cL}{\mathcal L}
\newcommand{\cM}{\mathcal M}
\newcommand{\cN}{\mathcal N}
\newcommand{\cO}{\mathcal O}

\newcommand{\cR}{\mathcal R}

\newcommand{\cV}{\mathcal V}

\newcommand{\GL}{\mathrm{GL}}
\newcommand{\HP}{\mathrm{H}}
\newcommand{\PH}{\mathrm{P}}

\renewcommand{\t}{\widetilde}

\newcommand{\tX}{\widetilde X}
\newcommand{\tY}{\widetilde Y}

\newcommand{\hbH}{\hat\bH}

\DeclareMathOperator{\Hom}{Hom}

\DeclareMathOperator{\Ker}{Ker}

\DeclareMathOperator{\spec}{Spec}
\DeclareMathOperator{\specan}{Specan}
\DeclareMathOperator{\proj}{Proj}
\DeclareMathOperator{\projan}{Projan}
\DeclareMathOperator{\pic}{Pic}

\DeclareMathOperator{\supp}{Supp}

\DeclareMathOperator{\Div}{Div}
\DeclareMathOperator{\LF}{LF}

\renewcommand{\:}{\colon}

\newcommand{\defset}[2]{{\left\{#1\,\left| \,#2 \right. \right\}}}
\newcommand{\X}{(X,o)}
\newcommand{\Y}{(Y,o)}

\newcommand{\ten}{\circle*{0.3}}
\newcommand{\bul}{\bullet}

\newcommand{\dis}{\displaystyle}

 \theoremstyle{plain}
\newtheorem{thm}{Theorem}[section]
\newtheorem{cor}[thm]{Corollary}
\newtheorem{lem}[thm]{Lemma}
\newtheorem{prop}[thm]{Proposition}

 \theoremstyle{definition}
\newtheorem{defn}[thm]{Definition}

\newtheorem{ex}[thm]{Example}

\newtheorem{ass}[thm]{Assumption}

\theoremstyle{remark}

\newtheorem{rem}[thm]{Remark}

\numberwithin{equation}{section}

\newcommand{\thmref}[1]{Theorem~\ref{#1}}
\newcommand{\lemref}[1]{Lemma~\ref{#1}}

\newcommand{\proref}[1]{Proposition~\ref{#1}}
\newcommand{\remref}[1]{Remark~\ref{#1}}

\newcommand{\defref}[1]{Definition~\ref{#1}}
\newcommand{\assref}[1]{Assumption~\ref{#1}}
\newcommand{\exref}[1]{Example~\ref{#1}}

\newcommand{\figref}[1]{Figure~\ref{#1}}
\newcommand{\sref}[1]{Section~\ref{#1}}

\begin{document}

\title{The geometric genus of splice-quotient singularities}

\author{Tomohiro Okuma}

\address{Department of Education, Yamagata University, 
 Yamagata 990-8560, Japan.}
  
\email{okuma@e.yamagata-u.ac.jp}

\thanks{Partly supported by the Grant-in-Aid for Young
Scientists (B), The Ministry of Education, Culture,
Sports, Science and Technology, Japan. 
}



\begin{abstract}
We prove a formula for  the geometric genus of
 splice-quotient singularities (in the sense of Neumann and 
 Wahl). This formula enables us to compute the invariant
 from  the resolution graph; in fact, it reduces the
 computation to that for splice-quotient singularities with
 smaller resolution graphs. 
We also discuss the dimension of the first cohomology groups of certain
 invertible sheaves on a  resolution of a  splice-quotient singularity.
\end{abstract}

\maketitle

\section{Introduction}

The topology (i.e., the link or the weighted dual graph) of
a normal surface singularity in general  
does not determine analytic invariants of the singularity.
It  is  challenging to see what kinds of analytic
invariants become  topological invariants under certain conditions, 
and to find a formula for computing the invariants from
resolution graphs.
The geometric genus $p_g$ is one of the fundamental analytic
invariants of singularities.
There are many studies on this invariant related to the
above issue; 
for example,  \cite{artin.rat},
\cite{la.me}, \cite{nw-casson}, \cite{nw-HSL},
\cite{nem.ellip}, \cite{swI}. 

Let  $X$ be a normal surface singularity whose link is a
$\Q$-homology sphere. Then there exists the ``universal abelian
cover'' $Y\to X$. 
Neumann and Wahl conjectured that  if $X$ is  $\Q$-Gorenstein,
then $Y$ is a complete intersection of ``splice type''
and  the geometric genus $p_g(X)$ is a topological
invariant.
Although this conjecture inspired research in the surface
singularity theory, counter-examples are now known (see
\cite{supiso}).
Splice type singularities are introduced by Neumann and
Wahl (\cite{nw-uac}, \cite{nw-CIuac}, \cite{nw-HSL}).
These singularities are a generalization of Brieskorn
complete intersections.
For a fixed resolution graph $\Gamma$, the associated splice
type singularities  form an equisingular family (see
\cite[Theorem 10.1]{nw-CIuac}, \cite[Theorem
4.3]{o.uac-certain}).
Therefore, $p_g$ of these singularities are the same and
determined by $\Gamma$ (this fact also follows from the result of
this paper). 
Neumann and Wahl have developed the theory of splice type
singularities, and recently proved the End-Curve Theorem
(\thmref{t:nw-end-curve}), 
which states that if a good resolution $\tX$ of $X$ satisfies the
``End-Curve condition'' (\defref{d:end-curve}), then $Y$ is
a splice type singularity associated with the resolution
graph of $X$ (the converse is also true). In this 
case, $X$ is  called a splice-quotient singularity.

Suppose that $X$ is a splice-quotient singularity with 
 resolution graph $\Gamma$, 
 and fix  a
component $E_v$ of the exceptional set $E \subset \tX$ such
 that $\delta:=(E-E_v)\cdot E_v\ge 3$. 
We consider a filtration $\{F_n\}_{n \ge 0}$ of $\cO_{X}$ associated with the
 prime divisor $E_v$, and singularities
 $\{X_i\}_{i=1}^{\delta}$ corresponding to connected
 components of the exceptional divisor $E-E_v$. Then $X_i$
 are also splice-quotient 
 singularities following  the End-Curve Theorem (\lemref{l:ind}). 
We can define an  invariant $c_v$ of the graded ring
 $\bigoplus _{n\ge 0}F_n/F_{n+1}$, which is sort of the constant term of
the  Hilbert polynomial; 
$c_v$ can be computed from $\Gamma$.
We prove the following.

\begin{thm}
If $X$ is a splice-quotient singularity, then 
$$
p_g(X)=c_v+\sum_{i=1}^{\delta}p_g(X_i).
$$
\end{thm}
By using this formula inductively, 
the computation of $p_g$ is reduced to that of
 $c_v$ (cf. \proref{p:last}).

We will also prove a formula for $h^1$ of invertible sheaves
on  $\tX$ related to the eigensheaves of $\cO_Y$ 
 (\thmref{t:genus}), which implies the formula for $p_g(Y)$.
In fact the above theorem is a corollary of this result.

 In \cite{CIC} A. N\'emethi and the author proved the Casson
 invariant conjecture of Neumann and Wahl (\cite{nw-casson},
 \cite[Theorem 6.3]{nw-HSL}) for splice type surface
 singularities
by applying our formula. 
 The conjecture can be reduced to proving an
 ``additivity property'' of the geometric genus under
 ``splicing'', while our formula relates an additivity
 property under ``plumbing.'' This gap is bridged by a new
 method in \cite[\S 4]{CIC}.

This paper is organized as follows.
  In \sref{s:pre}, we
 review basics of splice type singularities as  universal abelian
covers of normal surface singularities.
We show that a neighborhood of a connected 
  exceptional curve on  $\tX$ also satisfies the
  End-Curve condition.
This will enable us to use induction on the number of
  ``nodes'' of $E$.
In \sref{s:filtration}, we consider a weight 
filtration $\{I_n\}_{n\ge 0}$ of the local ring $\cO_{Y,o}$ 
with respect to  weights determined from
 the weighted dual graph of $E$;
this filtration is $\bH$-equivariant, where $\bH$ is the
 Galois group of the covering $q\: Y \to X$.
Then for any $\chi \in \Hom(\bH,\C^*)$,
the Hilbert series of the  $\chi$-eigenspace
  $\oplus I_n^{\chi}/I_{n+1}^{\chi}$ of the 
associated graded ring  is computed from $\Gamma$.
(The formula will be proved in Appendix.)  
By applying this fact and the Riemann-Roch formula, we  compute $h^1$ of certain invertible
 sheaves related to eigensheaves of $q_*\cO_Y$ (it also related to $I_n^{\chi}$'s).
In \sref{s:p_g}, we prove the main theorem.
By using the vanishing theorem,
 the desired invariant 
 is decomposed into one similar to $c_v$
  and invariants of
``smaller'' splice-quotient singularities.

Part of this paper was prepared during the author's stay at
the  R\'enyi Institute of Mathematics, Budapest, Hungary. 
The author would like to thank Professor Andr{\'a}s N\'emethi for
valuable discussions and comments, and the R\'enyi Institute for the
support and hospitality. 
He would also like to thank Professor Jonathan Wahl for his
helpful exposition  of the End-Curve Theorem, and the
referees for their careful reading and very helpful comments.

\section{Universal abelian covers and splice type singularities}
\label{s:pre}

We introduce some notations, and review some basics on
the splice type singularities and the universal abelian
covers of surface singularities.
Although a system of splice diagram equations, which defines a
splice type singularity, is  associated with  a weighted 
tree called a ``splice
diagram'' in origin (see \cite[\S 2]{nw-CIuac} for detail), 
 we construct them
in terms of ``monomial cycles'' on a resolution space of a surface
singularity (cf. \cite[\S 3]{o.uac-certain}, \cite[\S
13]{nw-CIuac}) for convenience of discussion. 

Let $\X$ be a germ of a normal complex surface singularity
and  $\pi\:\tX \to X$ a good resolution. 
Then the exceptional
divisor $E:=\pi^{-1}(o)$ has only simple normal crossings.
We denote  the link of the singularity $\X$ by $\Sigma$.
We may assume that $X$ is homeomorphic to a cone over
$\Sigma$.
In case $\pi$ is the minimal good resolution the weighted
dual graph of $\X$ is that of $E$.
It is known that the weighted dual graph of $\X$ and $\Sigma$ have the same 
information (\cite{neumann.plumbing}).

\begin{ass}\label{a:qhs}
We always assume that the link $\Sigma$ is a
 rational homology sphere;  it is 
 equivalent to that $E$ is a tree of rational curves.
In addition, we assume that $X$ is not a cyclic quotient
 singularity.
Therefore the weighted dual graph of  $E$ is not a chain (note that the universal abelian cover
 of a cyclic quotient  singularity is nonsingular).
\end{ass}

Let $\{E_v\}_{v \in \cV}$ denote the set of irreducible components
of $E$.  Let 
$$
\bL=\sum _{v\in \cV}\Z E_v \quad \text{and} \quad
\bL_{\Q}=\bL\otimes \Q.
$$
We call an element of $\bL$ (resp. $\bL_{\Q}$) a cycle
(resp. $\Q$-cycle).
Since the intersection matrix $(E_v\cdot E_w)$ is negative
definite, 
 for any $v \in \cV$ there exists a $\Q$-cycle
$E^*_v$ such that $E^*_v\cdot E_w=-\delta_{vw}$ for every $w
\in \cV$,  where
$\delta _{vw}$ denotes the Kronecker delta.
Let 
$$
\bL^*=\sum_{v \in \cV}\Z E_v^* \subset \bL_{\Q}.
$$
We may identify $\bL^*$ with $\Hom(\bL,\Z)$.
By the assumption,   $H_1(\Sigma, \Z)$ is a finite
group and we have a natural isomorphism
$$
 \bH:=\bL^*/\bL \to H_1(\Sigma, \Z)
$$ 
via  isomorphisms $\bL \cong H_2(\tX,
\Z)$ and $\bL^* \cong H_2(\tX, \Sigma,\Z)$ (cf. \cite[\S
2]{nem.line} or 
\cite[\S 2]{o.uac-rat}).
Here, we identify $\bH$ with $H_1(\Sigma, \Z)$.
The intersection pairing on $\bL_{\Q}$ induces a natural
pairing
\begin{equation}\label{eq:pairing}
 \bH\times \bH \to \Q/\Z.
\end{equation}
Let $\Div(\tX)$ denote the group of divisors on $\tX$.
We define a homomorphism
$$
c_1\: \Div(\tX) \to \bL^* \quad \text{by} \quad 
 c_1(D)=-\sum_{v}(D\cdot E_v)E^*_v.
$$
Clearly $c_1$ is surjective and $D\cdot E_v=c_1(D)\cdot E_v$
for all $v \in \cV$.

\subsection{The universal abelian cover of $X$}\label{ss:uac}

 There uniquely exists a finite morphism $q\:\Y \to \X$ of
 singularities that induces an  unramified Galois covering   
$Y \setminus \{o\} \to X    \setminus \{o\}$ with 
Galois group $\bH$. 
The  morphism $q\:X \to Y$ is called the {\itshape universal abelian
 covering} of $X$. 

In this subsection, we give an expression of the decomposition of $q_*\cO_{Y}$
into the $\bH$-eigensheaves, and a similar expression for the
structure sheaf of a partial resolution of $Y$; see
\cite{o.uac-rat} or \cite{nem.line} for details.
Let $\hat\bH=\Hom(\bH, \C^*)$. 
We define a homomorphism 
$$
\theta\: \bH \to \hbH \quad \text{by} \quad
\theta(h)(h')=\exp(2\pi\sqrt{-1}h\cdot h'),
$$
where $h\cdot h' \in \Q/\Z$ is determined by \eqref{eq:pairing}.
We also use the following notation:
\begin{itemize}
 \item $\theta(h,h'):=\theta(h)(h')$; then $\theta(\; ,
       \;)$ is symmetric.
 \item By abuse of notation, let $\theta$ also denote the composite $\bL^*
	\to \bH \xrightarrow{\theta} \hbH$.
 \item For $D \in \bL^*$ and $h \in \bH$, let $\theta(h,D):=\theta(D,h):=\theta(D)(h)$.
\end{itemize}
If $\Q$-divisors $D_1$ and $D_2$ are $\Q$-linearly
 equivalent, i.e., $nD_1\sim nD_2$ for some $n \in \N$, we
 write $D_1 \sim_{\Q}D_2$.
Note that if two integral divisors on $\tX$ are $\Q$-linearly
 equivalent then they are also linearly equivalent because
 $\pic (\tX)$ has no torsion by \assref{a:qhs}.
 
 There exists a set  $\{L_{\chi}\}_{\chi \in \hat\bH}$ of
 divisors on $\tX$ such that 
 \begin{enumerate}
  \item  $\theta(c_1(L_{\chi}))=\chi$,
  \item  $L_{\chi}\sim_{\Q} c_1(L_{\chi})$, and
  \item $[c_1(L_{\chi})]=0$, where $[ D ]$ denotes the
	integral part of a $\Q$-divisor $D$.
 \end{enumerate} 
Clearly, such an $L_{\chi}$ is uniquely determined up to
linear equivalence. 
Note that if $D \in \bL^*$ and $\chi=\theta(D)$,
 then $c_1(L_{\chi})-D \in \bL$.

We define a homomorphism
$$
\sigma\: \bL^* \to \Div(\tX) \quad \text{by} \quad 
\sigma(D)=L_{\theta(D)}+[D].
$$
Then $D\sim_{\Q}\sigma(D)$ for any $D \in \bL^*$. 
Clearly, $\sigma$ is a section of the homomorphism $c_1$.

Let $\cL_{\chi}:=\cO_{\tX}(-L_{\chi})$.

\begin{prop}[{cf. \cite[\S 3.2]{o.uac-rat}}]\label{p:uac-rat}
We have a collection $\{\cL_{\chi} \otimes \cL_{\chi'} \to
 \cL_{\chi\chi'}\}$ of homomorphisms that defines an 
$\cO_{\tX}$-algebra structure of
an $\cO_{\tX}$-module $\cA:=\bigoplus _{\chi \in 
 \hat\bH}\cL_{\chi}$ such that the following are satisfied.
\begin{enumerate}
 \item The projection $\specan _{X}\pi_*\cA \to X$
       coincides with $q\:Y \to X$.
 \item $\tY:=\specan _{\tX}\cA$  has only cyclic quotient
       singularity and there is a morphism $\rho\: \tY \to
       Y$, which is a partial resolution, and the following
       diagram is commutative:
$$
\begin{CD}
 \tY @>{p}>> \tX \\
@V{\rho}VV     @VV{\pi}V \\
Y @>>{q}> X
\end{CD}
$$
where $p$ is the natural projection.  Then, $p$ is unramified over
       $\tX\setminus E$. 
 \item The module $\cL_{\chi}$
       (resp. $\pi_*\cL_{\chi}$) is the
       $\chi$-eigensheaf of $p_*\cO_{\tY}$ (resp. $q_*\cO_{Y}$).
\end{enumerate}
\end{prop}

\subsection{Splice diagram equations}

Let $\delta_v=(E-E_v)\cdot E_v$ be the number of irreducible
components of $E$ intersecting $E_v$. A curve $E_v$ is called
an {\itshape end} (resp. a {\itshape node}) if $\delta_v=1$
(resp. $\delta_v\ge 3$). 
Let $\cE$ (resp. $\cN$) denote the set of indices of ends
(resp. nodes).
A connected component of $E-E_v$ is called a {\itshape branch} of $E_v$.

\begin{defn}
An element of a semigroup $\sum _{w \in \cE}\Z_{\ge 0}E^*_w$,
 where $\Z_{\ge 0}$ is the set of nonnegative integers, is
 called a {\itshape monomial cycle}.
Let $\C[z]:=\C[z_w; w \in \cE]$ be the polynomial ring in
 $\# \cE$  variables.
For a monomial cycle $D=\sum _{w \in \cE}\alpha_wE_w^*$, we
 associate a monomial  $z(D):=\prod_{w\in \cE} z_w^{\alpha_w} \in
 \C[z]$.
\end{defn}

\begin{defn}[Monomial Condition]
\label{c:A}
We say that $E$ (or its weighted dual graph) 
satisfies the monomial condition if for  any
 branch $C$ of any node $E_v$, 
 there exists a monomial cycle $D$
 such that $D-E^*_v$ is an effective integral cycle
 supported on $C$ (cf. \exref{ex:MC}).
In this case, $z(D)$ is called an {\itshape admissible monomial} belonging to the branch $C$.
\end{defn}

\begin{rem}
 The monomial condition is equivalent to the semigroup and
 congruence conditions (see \cite[\S
 13]{nw-CIuac}), 
which are required  for obtaining ``appropriate'' splice
 diagram equations (cf. \thmref{t:NW-CIuac}).  
Note that the original definition of admissible monomials
 requires only the semigroup condition (\cite{nw-CIuac},
 \cite{nw-HSL}).  
\end{rem}

\begin{defn}
 Assume that the monomial condition is satisfied.
Let $E_v$ be a node and let $C_1, \dots , C_{\delta_v}$ be
 the branches of $E_v$. 
Suppose   $\{m_1, \ldots ,m_{\delta_v}\}$  is a set of
 admissible monomials such that
 $m_i$ belongs to $C_i$ for  $i=1, \ldots ,\delta_v$. 
Let $F=(c_{ij})$, $c_{ij} \in \C$, 
be any $((\delta_v-2) \times  \delta_v)$-matrix such  
 that every maximal minor of it has rank $\delta_v-2$.
We define polynomials $f_1, \ldots ,f_{\delta_v-2}$ by 
$$
\begin{pmatrix}
 f_1\\ \vdots \\ f_{\delta_v-2}
\end{pmatrix}
= F
\begin{pmatrix}
 m_1\\ \vdots \\ m_{\delta_v}
\end{pmatrix}.
$$
We call the set $\{f_1 , \dots ,f_{\delta_v-2}\}$ a {\itshape Neumann-Wahl
 system}  at  $E_v$. 
Suppose that we have a Neumann-Wahl system $\cF_v$ at every
 node $E_v$.
Then we call the set 
$\cF:= \bigcup _{v \in \cN}\cF_v$ a {\itshape Neumann-Wahl system} associated
 with $E$.
Note that $\# \cF=\#\cE-2$.
\end{defn}

\begin{ex}
 \label{ex:MC}
Suppose that $E_4^2=-3$ and $E_w^2=-2$ for $w\ne 4$, and the
 weighted dual graph of $E$ is represented 
 as follows:
 \begin{center}
 \setlength{\unitlength}{0.5cm}
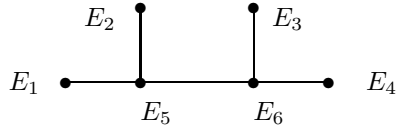
 \begin{picture}(11,3.5)(0,-1)
 \multiput(2,0)(2,0){2}{\ten}
 \multiput(7,0)(2,0){2}{\ten}
 \put(0.5,-0.2){$E_1$}
 \put(10,-0.2){$E_4$}
 \put(4,-1){$E_5$}
 \put(4,0){\line(0,1){2}} 
 \put(4,2){\ten}
 \put(7,0){\line(0,1){2}} 
 \put(7,2){\ten}
 \put(2.5,1.5){$E_2$}
 \put(7,-1){$E_6$}
 \put(7.5,1.5){$E_3$}
 \put(2,0){\line(1,0){7}} 
 \end{picture}
 \end{center}

Then the following equations show that  the monomial condition is
 satisfied. 
\begin{align*}
&2 E_1^*=E_5^*+E_1, \quad 2 E_2^*=E_5^*+E_2, \quad
 E_3^*+2E_4^*=E_5^*+(E_3+E_4+E_6), \\
&2 E_3^*=E_6^*+E_3, \quad 3 E_4^*=E_6^*+E_4,  \quad
E_1^*+E_2^*=E_6^*+(E_1+E_2+E_5).
\end{align*}
The corresponding admissible monomials form a Neumann-Wahl
 system 
$$
z_1^2+z_2^2+z_3z_4^2, \quad z_3^2+z_4^3+z_1z_2.
$$
\end{ex}

Let $(a_{vw})=-I^{-1}$, where $I$ denotes the intersection
matrix $(E_v\cdot E_w)$. 
Then every $a_{vw}$ is a positive rational
number and $E^*_v=\sum_{w \in \cV} a_{vw}E_w$.
We define positive integers $e_v$, $\ell_{vw}$, and $m_{vw}$ as follows:
$$
\ell_{vw}=\lvert \det I \rvert a_{vw}, \quad 
e_{v}=\lvert \det I \rvert /\gcd\defset{\ell_{vw}}{w \in
\cV}, \quad \text{and} \quad 
m_{vw}=e_va_{vw}.
$$
It is easy to see that $\gcd\{m_{vw}\}_{w \in \cV}=1$
for every $v \in \cV$.

\begin{defn}
 For any $v \in \cV$, we define the $v$-weight of the variable
 $z_w$, $w \in \cE$, to be $m_{vw}$.
Therefore, the  $v$-degree of a monomial
 $\prod_{w \in
 \cE}z_w^{\alpha_w}$ is $\sum_{w \in \cE}\alpha_wm_{vw}$.
Note that if $D=\sum _{w \in \cV}\beta_wE_w$ is a monomial cycle, then
 the  $v$-degree of $z(D)$ is equal to $e_v\beta_v=-e_vD\cdot E_v^*$.

Let $\C\{z\}:=\C\{z_w ; w \in \cE\}$ be the convergent power series
ring.  
Let $f=f_0+f_1 \in \C\{z\}$, where $f_0$ is a nonzero
 quasihomogeneous polynomial with respect to  
 the $v$-weight and  $f_1$ is a series in monomials of
 higher $v$-degrees. 
Then we call $f_0$ the {\itshape $v$-leading form} of $f$, and 
 denote it by $\LF_{v}(f)$.
We also define the {\itshape $v$-order} of $f$ to be the
 $v$-degree of $f_0$.
\end{defn}

\begin{defn}
We consider a finite set 
$$
 \defset{f_{vj_v}}{v \in \cN, \; j_v=1, \dots
 ,\delta_v-2}\subset \C\{z\}.
$$
If the set 
$$
 \defset{\LF_v(f_{vj_v})}{v \in \cN, \; j_v=1, \dots
 ,\delta_v-2}
$$ is a Neumann-Wahl system associated with $E$,  then  a 
 system of  equations  
$$
f_{vj_v}=0, \quad v \in \cN, \quad j_v=1, \dots ,\delta_v-2,
$$
is called the {\itshape splice diagram equations}.
A germ of a singularity defined by the splice diagram
 equations in $(\C^{\#\cE},o)$ is called a {\itshape splice type
 singularity}, or said to be of splice type.
\end{defn}

\begin{thm}[Neumann-Wahl {\cite[Theorem 2.6]{nw-CIuac}}]\label{t:from-nw-CIuac}
A splice type  singularity is  an isolated complete intersection surface
       singularity. 
\end{thm}

\subsection{$\bH$-action}\label{ss:H-action}

We define an action of $\bH$ on the power series ring
 $\C\{z\}$ as follows.
 For any monomial cycle $D$ and any element $h \in \bH$, we
 define $h \cdot z(D) \in \C \{z\}$ by
$$
h \cdot z(D)=\theta(h,D)z(D).
$$
This extends to an action on  $\C\{z\}$.
If $\cE=\{w_1, \dots ,w_n\}$,  the action corresponds to  a
 representation
$$
\bH \to \mathrm{U}(n), \quad h \to \mathrm{D}_h,
$$
where $\mathrm{D}_h$ denotes a diagonal matrix with $i$-th
 diagonal component $\theta(h,E_{w_i}^*)$.
If $\{f_{vj_v}\}$ is a Neumann-Wahl system associated with
 $E$, then for every $f_{vj_v}$ and $h \in \bH$,
we have $h\cdot f_{vj_v}=\theta(h,E_v^*)f_{vj_v}$, i.e.,
 $f_{vj_v}$ is in the $\theta(E_v^*)$-eigenspace of $\C\{z\}$.

\begin{thm}[Neumann-Wahl {\cite[Theorem 7.2]{nw-CIuac}}]\label{t:NW-CIuac}
 Suppose that $(Z,o)$ is a singularity defined by splice diagram equations
$$
f_{vj_v}=0, \quad v \in \cN, \quad j_v=1, \dots ,\delta_v-2,
$$
such that  $h\cdot  f_{vj_v}=\theta(h,E_v^*)f_{vj_v}$ for
 every $f_{vj_v}$ and $h \in \bH$.
Then we have the following.
\begin{enumerate}
\item $\bH$ acts freely on $Z\setminus \{o\}$, and thus,
      $Z':=Z/\bH$ is a normal surface singularity. 
\item The weighted dual graph of $Z'$ is the same as that of $X$.
\item The quotient map $Z\to Z'$ is the universal abelian covering.
\end{enumerate}
\end{thm}

\begin{defn}
 A singularity whose universal abelian cover is of splice
 type
(like $Z'$ in \thmref{t:NW-CIuac}) is
 called a {\itshape splice-quotient singularity}.
\end{defn}

\begin{defn}[End-Curve Condition]\label{d:end-curve}
We say that $\tX$ satisfies the end-curve condition if for
 each $w \in \cE$ there exists an irreducible curve $H_w
 \subset \tX$,
 not an exceptional curve, such that $H_w\cdot E=H_w\cdot E_w=1$
and $E_w^*+H_w\sim_{\Q} 0$ (this implies $e_w(E_w^*+H_w)\sim 0$).
In other words, 
 the end-curve condition is equivalent to that
$\cO_{\tX}(-\sigma(E_w^*))$  has no fixed component in $E$
 for every $w \in \cE$.  
In this case, a general section $s \in H^0(\cO_{\tX}(-\sigma(E_w^*)))$
 defines a divisor $\sigma(E_w^*)+H_w$, where $H_w$ is as
 above. 
We call such an $s$ an {\itshape end-curve section} of $E_w$.
\end{defn}

\begin{rem}
 If $\tX$ satisfies the end-curve condition, then so do the
 minimal good resolution and a resolution obtained by
 blowing up  $\tX$ at singular points of $E$ or a point
 $E_w\cap H_w$, $w \in \cE$.
\end{rem}

The following theorem is a generalization of \cite[Theorem
4.1]{nw-HSL}, which was announced in 
\cite[Theorem 6.1]{neumann.graph-3-mfd} and its proof will appear in
 {\cite{nw-preprint}.

\begin{thm}[End-Curve Theorem]\label{t:nw-end-curve}
 If $\tX$ satisfies the end-curve condition, then $E$
 satisfies the monomial condition and  $Y$ is of
 splice type. In fact, if  
$$
 \psi\: \C\{z_w; w \in \cE\} \to \cO_{Y,o}
$$
 is a homomorphism of $\C$-algebras that maps each $z_w=z(E_w^*)$
to an end-curve section of $E_w$, then $\psi$ is
 surjective and   $\Ker \psi$ is generated by functions
 $\{f_{vj_v}\}$
 as in \thmref{t:NW-CIuac}.
\end{thm}

The following proposition implies that if $X$ is a
splice-quotient singularity then any singularity obtained by
contracting a connected exceptional curve on $\tX$ is also a
splice-quotient singularity.

\begin{prop}\label{p:end-curve-induction}
 Let $v_0 \in \cE$, and let $\tX' \subset \tX$ be a sufficiently small
 neighborhood of the divisor $E':=E-E_{v_0}$.
If $\tX$ satisfies the end-curve condition, then so does $\tX'$.
\end{prop}
\begin{proof}
 Let $\cV'=\cV\setminus \{v_0\}$. For each $v \in \cV'$,
 define a $\Q$-cycle $E_v^{\times} \in \sum_{w \in
 \cV'}\Q E_w$ by the condition that $E_v^{\times}\cdot
 E_w=-\delta_{vw}$  for every $w \in \cV'$.
Assume that $E_{v_1}$ intersects $E_{v_0}$.
Since $E_{v_0}^{*}+H_{v_0}\sim_{\Q} 0$ for some curve $H_{v_0}$
 as in \defref{d:end-curve}, we have
 $E_{v_0}^{*}|_{\tX'} \sim _{\Q}0$.
We can easily see that
 $E_{v_0}^{*}-a_{v_0v_0}E_{v_0}=a_{v_0v_0}E_{v_1}^{\times}$. 
If $E_{v_1}$ is an end of $E'$, then we can take
 $E_{v_0}|_{\tX'}$ as  ``curve $H_{v_1}$.''
Suppose $w \in \cE\setminus\{v_0\}$.
Then $E_w^{*}+H_w\sim_{\Q}0$.
Since
 $a_{v_0v_0}E_w^*-a_{v_0w}E_{v_0}^*=a_{v_0v_0}E_w^{\times}$, 
we obtain that 
$$
 a_{v_0v_0}(E_{w}^{\times}+H_w|_{\tX'})=a_{v_0v_0}(E_{w}^{*}+H_w)|_{\tX'}
-a_{v_0w}E_{v_0}^*|_{\tX'}\sim_{\Q}0.
$$
Therefore, the end-curve condition is satisfied.
\end{proof}

\section{Filtration associated to a node}\label{s:filtration}

We use the notation of the preceding section.
Assume that the end-curve condition is satisfied and that
the universal abelian cover $p\: Y\to X$ is expressed  as in
\thmref{t:NW-CIuac} and \ref{t:nw-end-curve}.  
Note that  every monomial cycle $D$ determines a homogeneous
element  $\psi(z(D)) \in H^0(\cL_{\theta(D)})$ of
$\hbH$-graded local ring $\cO_{Y,o}$.

Throughout this section, we fix a node $E_v$. 
Let $C_1, \dots ,C_{\delta_v}$ denote the branches of $E_v$.

\begin{defn}
 For each $n \in \Z_{\ge 0}$,  $I_n$ denotes the
 ideal of $\cO_{Y,o}$ generated by the images of the
 elements of $\C\{z\}$ having  $v$-order $\ge n$.
Let  $\cG$ denote the associated graded algebra $\bigoplus 
_{n\ge 0}I_n/I_{n+1}$. Let $\cG_n=I_n/I_{n+1}$.
\end{defn}

\begin{thm}[Neumann-Wahl {\cite[Theorem 2.6]{nw-CIuac}}]\label{t:graded}
 Let $\{f_{wj_w}\}$ be the set
 of power series defining $Y$ as in \thmref{t:NW-CIuac}, and 
let  $I$ be the ideal of the polynomial ring $\C[z]$
 generated by the $v$-leading forms $\{\LF_v(f_{wj_w})\}$.
Then $\cG\cong \C[z]/ I$ and it is a reduced complete
 intersection ring. 
\end{thm}

\subsection{A geometric description of the filtration}

Let us recall the commutative diagram in \proref{p:uac-rat}.
Let $F:=p^{-1}(E)$.
Then $F$ is the $\rho$-exceptional set on $\tY$.
Let $F_v=p^{-1}(E_v)$ (this may be not
irreducible). Let $\pi'\: \tX \to X'$ (resp. $\rho'\: \tY
\to Y'$) be the 
morphism that contracts the divisor $E-E_v$
(resp. $F-F_v$) to normal points.
Then the natural morphism $p'\: Y' \to X'$ is  finite.
We have the  following commutative diagram:
$$
\xymatrix{
\tY \ar[rrr]^p \ar[dr]^{\rho'} \ar[dd]_{\rho} &&& \tX
 \ar[ld]_{\pi'} \ar[dd]^{\pi}\\
& Y' \ar[r]^{p'} \ar[dl]^{\rho_1}& X' \ar[rd] _{\pi_1} & \\
Y \ar[rrr]_q & && X
}
$$
where $\pi_1$  and $\rho_1$ are the natural morphisms.
Clearly, the exceptional sets of $\pi_1$  and $\rho_1$ are $E':=\pi'(E_v)$ and
$F':=\rho'(F_v)$, respectively. 
Since $p^*E_v=e_vF_v$ (see  \cite[Theorem 3.4]{o.uac-rat}),
it follows from the definition of the $v$-weight that
the ideal $I_n$ is generated by the images of power series
in monomials   $z(D)$ satisfying  $-e_vD\cdot E_v^*\ge n$
 or equivalently $p^*D\ge nF_v$.

We define a positive integer $a_v$ by 
$$
a_v:=e_v^2a_{vv}=e_vm_{vv}.
$$

\begin{lem}\label{l:filt1}
Let $n \in \Z_{\ge 0}$. Then, we have the following.
\begin{enumerate}
 \item $\cO_{\tX}(-\sigma(E_v^*))$ is $\pi$-generated and 
$\cO_{\tX}(-ne_vE_v^*)$ is generated by functions in $I_{na_v}\cap \cO_{X,o}$ near $E$.
 \item $I_n=\left(\rho_*\cO_{\tY}(-nF_v)\right)_o
=\left(\rho_{1*}\cO_{Y'}(-nF')\right)_o
=\left(\rho_*\cO_{\tY}([-n\rho'^*F'])\right)_o$, where
       $(\;)_o$  means the stalk at $o \in Y$.     
 \item $-F'$ is $\rho_1$-ample and $\rho_1$ coincides with
       the filtered blowing up 
$$
 \projan_Y\left(\bigoplus_{n\ge
       0}\rho_{1*}\cO_{Y'}(-nF')\right) \to Y.
$$
\end{enumerate}
\end{lem}
\begin{proof}
 Suppose $D_1$ and $D_2$ are monomial cycles such that each
 $z(D_i)$ is an admissible monomial belonging to the branch 
 $C_i$.
First, we note that $s_i:=\psi(z(D_i)) \in
 H^0\left(\cO_{\tX}(-\sigma(D_i))\right)$  has no zero
 outside $\bigcup_{D_i\cdot E_w\neq 0}E_w$ (see
 \defref{d:end-curve}).
Since  $\sigma(E_v^*) \le \sigma(D_i)$ and 
$\sigma(E_v^*) = \sigma(D_i)$ outside $C_i$ 
for $i=1,2$, the sections $s_1$ and $s_2$ 
 generate  $\cO_{\tX}(-\sigma(E_v^*))$.
We have the equality 
$$
 \cO_{\tX}(-\sigma(E_v^*))^{ne_v}=\cO_{\tX}(-ne_vE_v^*)\subset
 \cO_X
$$
in the
 algebra $p_*\cO_{\tY}$.  
As above, we see that sections  $s_1^{ne_v}$ and $s_2^{ne_v}$ generate
 $\cO_{\tX}(-ne_vE_v^*)$. 
Clearly, $ne_vD_i \in \bL$ and $-ne_vD_i\cdot E_v^*=na_{v}$. 
 Hence (1) follows.

The second equality in (2) and inclusion $I_{n}\subset
 \left(\rho_*\cO_{\tY}(-nF_v)\right)_o$ are clear.
The last equality follows from
the formula
 $\cO_{Y'}(-nF')=\rho'_*\cO_{\tY}([-n\rho'^*F'])$ (see
 \cite[(2.1)]{sakai.weil}).  
Since $\cG$ is reduced by \thmref{t:graded}, the argument of
 \cite[(2,2)]{tki-w} applies to our filtration; there exist
 valuations $V_1, \dots , V_t$ of the quotient field of $\cO_{Y,o}$ and
 rational numbers $q_1, \dots , q_t$ such that 
$$
 I_n=\defset{f \in \cO_{Y,o}}{V_i(f)\ge nq_i \text{ for  }
 1 \le i \le t} \quad \text{for all $n$. }
$$
Thus it suffices to show that there exists a positive
 integer $d$ such that
 $I_{md}=\left(\rho_{1*}\cO_{Y'}(-mdF')\right)_o$ for all $m$.
Let $F_v^*=\rho'^*F'$. Then we have an equality
 $a_vF_v^*=p^*(e_vE_v^*)$ of Cartier divisors.
It follows from  (1) that
 $\cO_{\tY}(-ma_vF_v^*)=I_{ma_v}\cO_{\tY}$, and thus,
$$
I_{ma_v}=\left(\rho_*\cO_{\tY}(-ma_vF_v^*)\right)_o
=\left(\rho_{1*}\cO_{Y'}(-ma_vF')  \right)_o  
$$
 for all $m$. 
This proves (2).

Since $\cO_{\tY}(-a_vF_v^*)$ is $\rho$-generated and
 trivial near $F-F_v$ but positive on $F_v$, the morphism
 $\rho_1$ is obtained by  blowing 
 up with respect to the ideal sheaf $\rho_*\cO_{\tY}(-ma_vF_v^*)$
 for some $m$.
Hence (3) follows.
\end{proof}

In the next lemma, we use the $a$-invariant of graded
rings induced by  Goto and Watanabe \cite[(3.1.4)]{G-W} (see
also \cite[3.6.13]{c-m}). 
By \thmref{t:graded} and formulas \cite[3.6.14--15]{c-m}, the
$a$-invariant $a(\cG)$ of $\cG$ is expressed as
$$
\sum_{w\in
 \cN}(\delta_w-2)m_{vw}-\sum_{w \in \cE}m_{vw}=\sum_{w\in
 \cV}(\delta_w-2)m_{vw}.
$$
By the definition of the $a$-invariant, the $n$-th graded component of $H^2_{\cG_+}(\cG)$
vanishes for $n>a(\cG)$, where $\cG_+=\bigoplus_{k>0}\cG_k$.

\begin{lem}\label{l:vanish}
 $H^1(\cO_{Y'}(-nF'))=0$ for $n > a(\cG)$.
\end{lem}
\begin{proof}
We apply the arguments of  Tomari-Watanabe \cite[\S 1]{tki-w}.
Let $\cR$ denote the Rees algebra $\bigoplus _{n\ge 0}I_nT^n\subset
 \cO_{Y,o}[T]$. 
 We may assume that $Y=\spec \cO_{Y,o}$ and $Y'=\proj \cR$.
Let $\cO_{Y'}(n)=\t{ \cR(n)}$.
Since  $\cO_{Y'}(n)$ is a divisorial sheaf for every $n\ge 0$   by Proposition 1.6
of \cite{tki-w},  it follows from \lemref{l:filt1} that 
 $\cO_{Y'}(n)\cong \cO_{Y'}(-nF')$.
On the other hand, (1.13) (ii) and (1.18) (i) of \cite{tki-w} imply that
$H^1(\cO_{Y'}(n))=0$ for $n>a(\cG)$.
\end{proof}

\subsection{The $\chi$-eigenspace of the filtration}

For any module $\cM$ with $\bH$-action and for any $\chi\in
\hbH$, let $\cM^{\chi}$ denote the $\chi$-eigenspace of
$\cM$ (see \S\ref{s:ap} for the definition).
The $\bH$-action on $\cO_{Y,o}$ induces an action on $I_n$.
For $\chi \in \hbH$, we have $I_n^{\chi}=I_n\cap
 \cO_{Y,o}^{\chi}$ and
 $\cG^{\chi}=I_n^{\chi}/I_{n+1}^{\chi}$. 

\begin{lem}\label{l:filt2}
 For any $D\in \bL^*$, the decomposition of the $\cO_{\tX}$-module
 $p_*\cO_{\tY}(-p^*D)$ into eigensheaves is given by
$$
p_*\cO_{\tY}(-p^*D)
=\bigoplus_{\chi\in\hbH}\cO_{\tX}(-L_{\chi}+[c_1(L_{\chi})-D])
$$
\end{lem}
\begin{proof}
 It follows from Lemma 3.2 and the proof of Theorem 3.4 of \cite{o.uac-rat}.
\end{proof}

We define  invertible sheaves  $\cL_{\chi,n}$ and
$\cL_{\chi,n}^*$ by
\begin{align*}
 \cL_{\chi,n}&=\cO_{\tX}(-L_{\chi}+[c_1(L_{\chi})-(n/e_v)E_v]), \\
 \cL_{\chi,n}^*&=\cO_{\tX}(-L_{\chi}+[c_1(L_{\chi})-(n/m_{vv})E_v^*]).
\end{align*}

\begin{lem}\label{l:filt3}
For any $\chi \in \hbH$ and $n \in \Z_{\ge 0}$, we have the
 following.
 \begin{enumerate}
  \item There exists a unique minimal effective cycle
	$D_{\chi,n} \in \bL$ such that
$$
\cL_{\chi,n}(-D_{\chi,n}):=\cL_{\chi,n}\otimes\cO_{\tX}(-D_{\chi,n})
$$	
is $\pi$-nef.
 \item 
$\dis	I_n^{\chi}=\left(\pi_*\cL_{\chi,n}\right)_o
=\left(\pi_*\cL_{\chi,n}^*\right)_o  
=\left(\pi_*\cL_{\chi,n}(-D_{\chi,n})\right)_o.$
 \end{enumerate}
\end{lem}
\begin{proof}
(1) and the  equality $\pi_*\cL_{\chi,n}=\pi_*\cL_{\chi,n}(-D_{\chi,n})$
 follows from  \cite[4.2]{nem.line}.
Note that $\rho'^*F'=m_{vv}^{-1}p^*E_v^*$.
By \lemref{l:filt1} and
 \ref{l:filt2}, we obtain other equalities.
\end{proof}

\begin{rem}
Suppose  $c_1(e_vL_{\chi})=\sum _{w \in\cV}k_wE_w$.
 Then  $k_v$ is an integer such
 that $0 \le k_v<e_v$.
By (2) of \lemref{l:filt3}, we have
$I_n=I_{n+1}$ if $n\not\equiv k_v \pmod{e_v}$; therefore,
$\cG^{\chi}=\bigoplus_{n \ge  0}\cG_{k_v+e_vn}^{\chi}$. 

\end{rem}

For any $\chi \in \hbH$, let $\HP^{\chi}(t)$ denote the
Hilbert series of the 
graded module $\cG^{\chi}$, i.e.,
$$
\HP^{\chi}(t)=\sum_{i\ge 0}(\dim \cG^{\chi}_i)t^i.
$$
We define a function $\PH^{\chi}$ on $\N$ by
$$
\PH^{\chi}(n)=\sum_{i=0}^{n-1}(\dim \cG^{\chi}_i).
$$

\begin{prop}[see \S\ref{s:ap}. Appendix]\label{p:hp}
We have the following formula
$$
 \HP^{\chi}(t)=\frac{1}{|\bH|}\sum _{h\in  \bH}
\chi^{-1}(h) \prod_{w\in \cV}
\left(1-\theta(h, E^*_w)t^{m_{vw}}\right)^{\delta_w-2}.
$$
In particular,  $\HP^{\chi}(t)$ and $\PH^{\chi}(n)$ are computed from the
weighted dual graph of $E$.
\end{prop}

For any cycle $D$ and an invertible sheaf
 $\cL$ on $\tX$, $\chi(\cL\otimes \cO_D)$ denotes the Euler
 characteristic, i.e., $h^0(\cL\otimes \cO_D)-h^1(\cL\otimes
 \cO_D)$. 
By the Riemann-Roch formula, 
$$
\chi(\cL\otimes \cO_D)=-D\cdot (D+K_{\tX})/2+\cL\cdot D.
$$

\begin{thm}\label{t:h1eigensheaves}
 For any $\chi \in \hbH$ and $n \in \Z_{\ge 0}$, and any
 effective cycle $D\le D_{\chi,n}$, we have the
 following:
\begin{gather*}
 \dim H^0(\cL_{\chi})/H^0(\cL_{\chi,n}(-D))=\PH^{\chi}(n),
 \quad \text{and} \\
h^1(\cL_{\chi,n}(-D))=\chi\left( \cL_{\chi} \otimes \cO_{D'}
 \right)-\PH^{\chi}(n) +h^1(\cL_{\chi}),
\end{gather*}
where $D'=D-[c_1(L_{\chi})-(n/e_v)E_v]$.
Furthermore, these are computed from the weighted dual graph
 of $E$ and the intersection numbers $\cL_{\chi}\cdot E_w$. 
\end{thm}
\begin{proof}
It follow from \lemref{l:filt3}  that
$$
 \dim H^0(\cL_{\chi})/H^0(\cL_{\chi,n}(-D))=\dim
 I_0^{\chi}/I_n^{\chi}=\PH^{\chi}(n). 
$$
We have $D'\ge 0$ since $[c_1(L_{\chi})]=0$.
From the exact sequence
$$
0 \to \cL_{\chi,n}(-D) \to \cL_{\chi} \to 
\cL_{\chi} \otimes  \cO_{D'} \to 0,
$$
we have the second formula.
Finally, we have to show $h^1(\cL_{\chi})$
 can be computed from the graph; however, it
 follows from \thmref{t:genus}.
\end{proof}

\section{The formulas}\label{s:p_g}

The situation is the same as in the preceding section.
Therefore, $\tX$ satisfies the end-curve condition and the node $E_v$ is again fixed.
Let  $\tX_i\subset \tX$ be a sufficiently small neighborhood
of the branch $C_i$ of $E_v$.
Suppose that the irreducible components of $C_i$ are indexed
by a set $\cV_i\subset \cV$. 
Let $E^*_{w,i}$ denote the $\Q$-cycle supported on $C_i$
 satisfying $E^*_{w,i}\cdot E_{w'}=-\delta _{ww'}$ for every
 $w'\in \cV_i$. 
We define groups $\bL_i$, $\bL^*_i$, and $\bH_i$ as follows:
$$
\bL_i=\sum_{w\in\cV_i}\Z E_w, \quad
\bL^*_i=\sum_{w\in\cV_i}\Z E^*_{w,i}, \quad \bH_i=\bL^*_i/\bL_i.
$$
A map $\theta_i\: \bH_i \to \hbH_i$ is defined as $\theta$ in
\S \ref{ss:uac}; we will follow the notational convention
used there.
For each $i$, we define a map
$$
\phi_i\: \bL^* \to \bL^*_i \quad \text{by} \quad 
\sum_{w \in  \cV} \alpha_wE^*_w \mapsto \sum_{w \in  \cV_i} \alpha_wE^*_{w,i}. 
$$
Recall that the natural map $\{c_1(L_{\chi})\}_{\chi \in
 \hbH}\to \bH$ is  bijective.
We define
$$
\psi_i\: \hbH \to \hbH_i \quad \text{by} \quad \chi \mapsto
 \theta_i(\phi_i(c_1(L_{\chi}))). 
$$

\begin{rem}\label{r:phi-psi}
 We have $L_{\chi}|_{\tX_i} \sim_{\Q}
 \phi_i(c_1(L_{\chi}))$ (cf. the proof of \proref{p:end-curve-induction}).
\end{rem}

There is a set $\{L_{\lambda}\}_{\lambda \in \hat\bH_i}$ of
 divisors on $\tX_i$ having properties similar to those of
 $\{L_{\chi}\}_{\chi \in \hat\bH}$ (see \S \ref{ss:uac}). 
For $\chi \in \hbH$, let $D_{\chi,i} =-[\phi_i(c_1(L_{\chi}))]$.
Then
\begin{equation}\label{eq:phi}
  L_{\psi_i(\chi)}\sim_{\Q}\phi_i(c_1(L_{\chi}))+D_{\chi,i}.
\end{equation}
In the following, if $\cL$ is a sheaf then $\chi(\cL)$
 denotes the Euler  characteristic of $\cL$.

\begin{lem}\label{l:induction}
For any  $\chi\in \hbH$ and $1\le i\le \delta_v$, we have the  following:
\begin{enumerate}
 \item $D_{\chi,i} \ge 0$ and
       $L_{\chi}|_{\tX_i}\sim L_{\psi_i(\chi)}-D_{\chi,i}$.   
 \item $H^0(\cO_{\tX_i}(-L_{\chi}))=H^0(\cO_{\tX_i}(-L_{\psi_i(\chi)}))$.
 \item       $h^1(\cO_{\tX_i}(-L_{\chi}))
=h^1(\cO_{\tX_i}(-L_{\psi_i(\chi)}))-\chi(\cO_{D_{\chi,i}}(-L_{\chi}))$. 
\end{enumerate}
\end{lem}
\begin{proof}
We write $L_i=L_{\psi_i(\chi)}$ and $D_i=D_{\chi,i}$.
Suppose that $E_{v_1}\subset C_i$ intersects $E_v$,
and write
 $c_1(L_{\chi})=\beta E_v+F_1+F_2$, where $\supp(F_1)\subset C_i$
 and $\supp(F_2)\subset \bigcup_{j\neq i}C_j$.
Then $\phi_i(c_1(L_{\chi}))=-\beta E_{v_1,i}^*+F_1$.
Since $[\beta E_v+F_1]=0$, we have 
$$
 D_i=-[-\beta E_{v_1,i}^*+F_1]\ge
 0.
$$
By \remref{r:phi-psi}  and \eqref{eq:phi}, we obtain  (1).
We have the natural inclusions
$$
 H^0(\cO_{\tX_i}(-L_i))\subset H^0(\cO_{\tX_i}(-L_i+D_i)) \subset
 H^0(\tX_i\setminus C_i, \cO_{\tX_i}(-L_i)).
$$
However, $H^0(\cO_{\tX_i}(-L_i))=H^0(\tX_i\setminus C_i, \cO_{\tX_i}(-L_i))$ by
 \cite[Proposition 4.3 (1)]{o.uac-rat}.
Therefore,  (2) follows from (1).
Since $D_i\ge 0$, we have the following exact sequence
$$
0 \to \cO_{\tX_i}(-L_i) \to \cO_{\tX_i}(-L_i+D_i) \to
 \cO_{D_i}(-L_i+D_i) \to 0.
$$
Then, (3) follows from  (1) and (2).
\end{proof}

We define Cartier divisors $C$, $C'$, $D$, and $D'$  as
follows:
$$
C:=e_vE_v^*, \quad C':=\pi' _*C, \quad D:=p^*C, \quad D':=\rho'_*D.
$$
Then $D'=a_v F'$, where $F'$ is the reduced exceptional
divisor on $X'$. 
Since $\rho'^*D'=D$, it follows from  \lemref{l:filt1} and
\ref{l:filt3}
that 
$$
I_{ma_v}=\left(\rho_*\cO_{\tY}(-mD)\right)_o \quad
 \text{and} \quad 
I_{ma_v}^{\chi}=\left(\pi_*\cO_{\tX}(-L_{\chi}-mC)\right)_o.
$$

\begin{lem}\label{l:vanish-eigen}
 If $m \in \N$ satisfies $m >a(\cG)/a_v$, 
then for every $\chi \in \hbH$,
$$
H^1(\pi'_*\cO_{\tX}(-L_{\chi}-mC))=0.
$$
\end{lem}
\begin{proof}
By the equalities $D=\rho'^*D'=p^*C$ and \lemref{l:filt2},
$$
p'_*\cO_{Y'}(-mD')=\pi'_*p_*\cO_{\tY}(-mD)\cong \bigoplus
 _{\chi \in \hbH} \pi'_*\cO_{\tX}(-L_{\chi}-mC).
$$
 Since $p'$ is finite, it follows from \lemref{l:vanish} that
$$
H^1(p'_*\cO_{Y'}(-mD'))=H^1(\cO_{Y'}(-ma_vF'))=0.
$$
\end{proof}

\begin{lem}\label{l:onbranches}
 Assume that $m \in \N$ satisfies $m >a(\cG)/a_v$. For every
 $\chi \in \hbH$, we have
$$
h^1(\cO_{\tX}(-L_{\chi}-mC))=\sum_{i=1}^{\delta_v}h^1(\cO_{\tX_i}(-L_{\chi})).
$$
\end{lem}
\begin{proof}
Let $\cL=\cO_{\tX}(-L_{\chi}-mC)$.
 From the spectral sequence 
$$
 E_2^{i,j} = H^i(R^j\pi'_*\cL) \Rightarrow  H^n(\cL),
$$ 
 we have the following exact sequence:
$$
0 \to H^1(\pi'_*\cL) \to
 H^1(\cL)
\to H^0(R^1\pi'_*\cL) \to 0.
$$
However, \lemref{l:vanish-eigen} implies that $H^1(\cL)\cong H^0(R^1\pi'_*\cL)$.
Let $x_i\in X'$ denote the singularity obtained by
 contracting $C_i$.
Then the support of $R^1\pi'_*\cL$ is in the set $\{x_i\}_i$.
Since $C=\pi'^*C'$ and $C'$ is a Cartier divisor,
$$
\left(R^1\pi'_*\cL\right)_{x_i}\cong \left(R^1\pi'_*\cO_{\tX}(-L_{\chi})
 \otimes\cO_{X'}(-mC')\right)_{x_i} \cong
 \left(R^1\pi'_*\cO_{\tX}(-L_{\chi})\right)_{x_i}. 
$$ 
This proves the lemma.
\end{proof}

Let $K$ denote the canonical divisor on $\tX$.

\begin{thm}\label{t:genus}
 For any  $\chi \in \hbH$ and any $m \in \N$ greater than
 $a(\cG)/a_v$, we have
\begin{multline*}
h^1(\cL_{\chi})=\PH^{\chi}(ma_v)-\frac{1}{2}(m^2a_v-me_v(K+2L_{\chi})\cdot
 E_v^*) \\
+\sum_{i=1}^{\delta_v}\left(h^1(\cO_{\tX_i}(-L_{\psi_i(\chi)}))
-\chi(\cO_{D_{\chi,i}}(- L_{\chi})\right). 
\end{multline*}
Furthermore, $h^1(\cL_{\chi})$ is a topological invariant;
 in fact, it is computed from the weighted
 dual graph $\Gamma$ of $E$ and the intersection numbers $L_{\chi}\cdot E_v$.
\end{thm}

The formula above directly induces a formula for $p_g(Y)$
because
$$
p_g(Y)=h^1(\cO_{\tY})=\sum _{\chi \in \hbH}h^1(\cL_{\chi}).
$$

\begin{cor}\label{c:p_g}
 $p_g(X)$ and $p_g(Y)$ can be computed from $\Gamma$. In
 fact, for  $m >a(\cG)/a_v$, we have 
$$
p_g(X)=\PH^{1}(ma_v)-\frac{1}{2}(m^2a_v-me_vK\cdot
 E_v^*) 
+\sum_{i=1}^{\delta_v}p_g(X',x_i). 
$$
\end{cor}

\begin{lem}\label{l:ind}
Every $\tX_i$ satisfies the end-curve condition.
\end{lem}
\begin{proof}
It follows from \proref{p:end-curve-induction}.
\end{proof}

\begin{proof}[Proof of \thmref{t:genus}]
 From the exact sequence 
$$
0 \to \cO_{\tX}(-L_{\chi}-mC) \to \cO_{\tX}(-L_{\chi}) \to
 \cO_{mC}(-L_{\chi}) \to 0
$$
we obtain 
$$
h^1(\cL_{\chi})=\dim
 \frac{H^0(\cO_{\tX}(-L_{\chi}))}{H^0(\cO_{\tX}(-L_{\chi}-mC))}
-\chi\left(\cO_{mC}(-L_{\chi})\right)+h^1(\cO_{\tX}(-L_{\chi}-mC)).
$$
Since $C=e_vE_v^*$ and $a_v=e_v^2a_{vv}$,  using
the Riemann-Roch formula, we have
\begin{align*}
 \chi\left(\cO_{mC}(-L_{\chi})\right) &
=-\frac{1}{2}mC\cdot(mC+K)-mL_{\chi}\cdot C \\
& =\frac{1}{2}m^2a_v-\frac{1}{2}me_v(K+2L_{\chi})\cdot
 E_v^*.
\end{align*}
On the other hand, by \lemref{l:filt3} we have
$$
\dim
 \frac{H^0(\cO_{\tX}(-L_{\chi}))}{H^0(\cO_{\tX}(-L_{\chi}-mC))}
=\dim I_0^{\chi}/I_{ma_v}^{\chi} =\PH^{\chi}(ma_v).
$$
The formula now follows from Lemmas \ref{l:induction} and \ref{l:onbranches}.

Next, we  show the second assertion by induction on $\#\cN$.
First, note that the invariant
 $h^1(\cL_{\chi})
-\sum_{i=1}^{\delta_v}h^1(\cO_{\tX_i}(-L_{\psi_i(\chi)}))$
 is computable from $\Gamma$.  
If $\#\cN=1$, then each $C_i$ is a chain.
Hence, every $h^1(\cO_{\tX_i}(-L_{\psi_i(\chi)}))$ is
 computable from $\Gamma$  (cf. \cite[4.3]{nem.line}).
For $\#\cN>1$, the invariant $h^1(\cO_{\tX_i}(-L_{\psi_i(\chi)}))$ is
 again computable from $\Gamma$ by \lemref{l:ind} and the hypothesis of the induction.
\end{proof}

In the formula of  \thmref{t:genus},
$$
c_v^{\chi}:= \PH^{\chi}(ma_v)-\frac{1}{2}(m^2a_v-me_v(K+2L_{\chi})\cdot
 E_v^*)$$ 
is independent of $m\gg 0$. 
In fact, $\PH^{\chi}(mka_v)$ is a polynomial function of $m$ for every integer $k>a(\cG)/a_v$.
If $\PH^{\chi}(ml)$ is also a polynomial function of $m$ for
some $l \in \N$, then its constant term coincides with $c_v^{\chi}$.
Thus, $c_v^{\chi}$ is an invariant of the series $\HP^{\chi}(t)$. 
\begin{prop}\label{p:last}
If $\HP^{\chi}(t)$ is expressed as  $p(t)+r(t)/q(t)$, 
where $p$, $q$, and $r$ are polynomials with $\deg r < \deg q$, then
 $c_v^{\chi}=p(1)$.
\end{prop}
\begin{proof}
We define a function $f$ on $\Z_{\ge 0}$ by $\sum_{n\ge 0} f(n) t^n=
 r(t)/q(t)$. For a positive integer $l$, let $P_l(m)=
 \sum_{n=0}^{ml-1}f(n)$. 
Then it is sufficient to prove that $P_l(m)$ is a polynomial
 function with $P_l(0)=0$ for some $l>\deg p$. 
We may assume that $q$ and $r$ have no common root.
By  \proref{p:hp}, there exists a positive integer $N>\deg p$
 such that $\alpha^N=1$ for every root $\alpha$ of $q(t)$.
Therefore it follows from \cite[Corollary 1.5]{stan.comb} that there exist
 polynomials $f_0, \dots ,f_{N-1}$ such that $f(n)=f_i(n)$
 if $n\equiv i \pmod N$. Then
 $P_N(m)=\sum_{i=0}^{N-1}\sum_{n=0}^{m-1}f_i(nN+i)$. Since
 $f_i(nN+i)$ is a polynomial  function of $n$ for each  $i$,
 $\sum_{n=0}^{m-1}f_i(nN+i)$ is a  polynomial function of
 $m$ whose constant term is zero. This proves the assertion.
\end{proof}

\begin{ex}
Let us consider a weighted dual graph $\Gamma$ represented as in
 \figref{fig:Gamma}.
\begin{figure}[htb]
 \begin{center}
 \setlength{\unitlength}{0.5cm}
 \begin{picture}(20,3.5)(-1,-1)
 \multiput(0,0)(2,0){10}{\ten}
 \put(-1,0.2){$w_1$}
 \multiput(-0.8,-1)(2,0){2}{$-2$}
 \put(3.2,-1){$-1$}
 \put(4.2,0.2){$v_1$}
 \put(4,0){\line(0,1){2}} 
 \put(4,2){\ten}
 \put(4.2,2){$-4$}
 \put(3,1.5){$w_2$}
 \put(5.2,-1){$-16$}
 \put(8.2,0.2){$v_0$}
 \put(7.2,-1){$-2$}
 \put(8,0){\line(0,1){2}} 
 \put(8,2){\ten}
 \put(8.2,2){$-2$}
 \put(7,1.5){$w_3$}
 \put(9.2,-1){$-4$}
 \multiput(11.2,-1)(2,0){4}{$-2$}
 \put(0,0){\line(1,0){18}} 
 \put(13.2,0.2){$v_2$}
 \put(16,1.5){\ten}
 \put(18,1.5){\ten}
 \put(14,0){\line(4,3){2}} 
 \put(16,1.5){\line(1,0){2}} 
 \put(15.8,2){$-2$}
 \put(17.8,2){$-2$}
 \put(18.5,1.2){$w_4$}
 \put(18.5,0){$w_5$}
 \end{picture}
 \caption{\label{fig:Gamma}}
 \end{center}
 \end{figure}
We can show that this graph  satisfies the
 monomial condition. 
Indeed, the set of the following functions, $z_i$ being the variable
 corresponding to the end $w_i$, is a Neumann-Wahl system
 associated with $\Gamma$:
$$
z_1+z_3^2+z_4z_5 \text{ (at $v_0$)}, \quad
 z_1^3+z_2^4+z_3^{30}  \text{ (at $v_1$)}, \quad 
z_1z_2^7z_3+z_4^3+z_5^3  \text{ (at $v_2$)}.
$$
(This is an example of a splice type singularity of
 embedding dimension  less than $\#\cE$.) 
Now, we suppose that $X$ is a splice-quotient singularity with  weighted
 dual graph  $\Gamma$.
Then, $X$ is numerical Gorenstein (i.e., $c_1(K_{\tX})\in \bL$) and 
 $p_a(Z)=1+Z\cdot (Z+K_{\tX})/2=4$, where $Z$ is
 Artin's fundamental cycle.
We can show that $X$ is actually a Gorenstein singularity, i.e.,
 $K_{\tX} \sim c_1(K_{\tX})$ holds. 
Since any splice-quotient singularity is $\Q$-Gorenstein, $K_{\tX} \sim_{\Q}
c_1(K_{\tX})$. On the other hand, the condition
 $H^1(\cO_E)=0$ implies that $\pic (\tX)$ has no torsion. 
Therefore  $K_{\tX} \sim c_1(K_{\tX})$.

Since  $|\bH|=36$, the universal abelian covering of $X$ is not
 trivial.
In the following we  compute the geometric genus of $X$
 by applying our formula (to do this we will need a computer
 algebra system).

Recall that the group $\bH$ is generated by $\{E_w^*\}_{w
 \in \cV}$ and
 the  relations on the generators are given by the
 intersection  matrix $I=(E_v\cdot E_w)$. By taking the Hermite normal
 form of $I$, we have a simplified expression of $\bH$; for example,
$\bH$ has a set of generators $\{E_{w_2}^*, E_{w_3}^*,E_{w_4}^*\}$
 with relations $2E_{w_2}^*, 6E_{w_3}, 
 E_{w_2}^*+3E_{w_3}^*+3E_{w_4}^* \in \bL$.
Let  $H_{v_0}(t)$ denote  the  Hilbert series associated with
 the node $v_0$. 
Let $\Lambda=\defset{\lambda\in (\Z_{\ge 0})^3}{\lambda_1<2,
 \lambda_2<6, \lambda_3<3}$ and
 $E^*_{\lambda}=\sum\lambda_iE^*_{w_{i+1}}$ for $\lambda\in \Lambda$. 
Then,
\begin{align*}
 H_{v_0}(t) &=\frac{1}{36}\sum_{\lambda\in \Lambda}\prod_{w\in
 \cV}\left(1-\exp (2\pi\sqrt{-1} E^*_{\lambda}\cdot
 E^*_w)t^{m_{v_0w}}\right)^{\delta_w-2} \\
&=\frac{t^{24}-t^{21}+t^{18}-t^{15}+3
   t^{12}-t^9+t^6-t^3+1}{t^{15}-t^{12}-t^3+1} \\
 &=\frac{4 t^{12}-2 t^9+2 t^6-2 t^3+1}{t^{15}-t^{12}-t^3+1}+t^9+t^3.
\end{align*}
By \proref{p:last}, we have $c_{v_0}=2$.
Next, let $\Gamma_i$ denote the branch of $v_0$ containing
 $v_i$ and $H_{\Gamma_i,v_i}(t)$  the  Hilbert series associated with
 the node $v_i$ on $\Gamma_i$ ($i=1,2$).
Then, we have
\begin{align*}
 H_{\Gamma_1,v_1}(t)&=\frac{t^{36}-t^{33}+t^{24}-t^{18}+t^{12}-t^3+1}
{t^{19}-t^{16}-t^3+1} \\
&=\frac{t^{12}+t^8+t^4-t^3-t^2-t+1}{t^{19}-t^{16}-t^3+1}+t^{17}+t^5+t^2+t
 \\
 H_{\Gamma_2,v_2}(t)&=\frac{t^{24}+1}{t^{20}-t^{14}-t^6+1}
=\frac{t^{18}+t^{10}-t^4+1}{t^{20}-t^{14}-t^6+1}+t^4.
\end{align*}
Therefore $c_{\Gamma_1,v_1}=4$ and $c_{\Gamma_2,v_2}=1$ (the
 second also follows from that $\Gamma_2$ corresponds to a
 minimally elliptic singularity).
Now, we conclude that 
$$
 p_g(X)=c_{v_0}+c_{\Gamma_1,v_1}+c_{\Gamma_2,v_2}=7.
$$
\end{ex}

\section{Appendix: Molien's formula}\label{s:ap}

Let $S=\C[z_1, \dots ,z_n]$ be the polynomial ring in $n$
variables. Suppose that the grading on $S$ is given by 
$\deg z_j=w_j\in \N$.
Let $S_i \subset S$ denote the $i$-th graded component.
Let $G$ be a finite subgroup of $\GL(n,\C)$ consisting of
diagonal matrices. 
Then  $G$ acts linearly on $S$ in the usual manner and
every graded component $S_i$ is invariant under the action.
We denote  the action by
$(g,f)\mapsto g\cdot f \in S$ for $(g,f) \in G\times S$. 
For any character $\chi \in \hat G=\Hom(G,\C^*)$, the
$\chi$-eigenspace $S^{\chi}$ is the space
$$
\defset{f \in S}{g\cdot f=\chi(g)f \text{ for all }
g \in G}.
$$
Then $S^{\chi}$ is also graded. Let $S^{\chi}_i=S^{\chi}\cap
S_i$. 
We denote by $\cH^{\chi}(t)$ the Hilbert series
of $S^{\chi}$, i.e.,
$$
\cH^{\chi}(t)=\sum_{i\ge 0}(\dim S^{\chi}_i) t^i.
$$
Then as in the  proof of \cite[Theorem 6.4.8]{c-m},  
we have the following

\begin{thm}[Molien's formula]\label{t:molien}
If $g_j$ denotes the $j$-th diagonal component of $g$, 
$$
\cH^{\chi}(t)=\frac{1}{|G|} \sum _{g \in
 G}\frac{\chi(g)^{-1}}{\prod_{j=1}^n(1-g_jt^{w_j})}.
$$ 
\end{thm}

Next, we consider a complete intersection case.
 Let $f_1, \dots ,f_m \in S$ be a regular sequence of homogeneous elements with
 $\deg (f_i)=d_i$, and $I:=(f_1, \dots ,f_m)\subset S$.
Assume that there exist $\chi_1, \dots ,\chi_m \in \hat G$ such that
 $g\cdot f_i=\chi_i(g)f_i$ for all $g \in G$.
Then $S/I$ is a graded $\C$-algebra with a natural
 $G$-action, and every
 graded  component $(S/I)_i$ is 
 invariant under this action. 
Let $H^{\chi}(t)$ denote the Hilbert series of  $\chi$-eigenspace
 $(S/I)^{\chi}$.

\begin{thm}\label{t:molien2}
We have the formula
$$
 H^{\chi}(t)=\frac{1}{|G|}\sum _{g \in G}
 \chi^{-1}(g)\frac{\prod_{i=1}^m (1-\chi_i(g)t^{d_i})}
{\prod_{j=1}^n(1-g_jt^{w_j})}. 
$$
\end{thm}

Let us outline the proof.
We consider the Koszul complex of the sequence
$f_1, \dots ,f_m$. 
We use the following notation.
\begin{itemize}
  \item Let $\Gamma_p=\defset{(i_1, \dots ,i_p)}{i_j \in \N,
	\; 1\le i_1< \dots <i_p\le m }$.
 \item For $\gamma=(i_1, \dots ,i_p) \in \Gamma_p$, we set 
	$d_{\gamma}:=\sum _{j=1}^pd_{i_j}$,
	$\chi_{\gamma}:=\prod_{j=1}^p \chi_{i_j}$ and 
	$\gamma\setminus i_j:=(i_1, \dots , \hat {i_j},\dots
	,i_p)$, where $\hat  {i_j}$ means that  $i_j$ is to be omitted.
 \item Let $K_p$ denote a graded  $S$-module with free basis
       $v_{\gamma}$, $\gamma  \in \Gamma_p$ defined as
 $$
 K_p=\bigoplus_{\gamma \in \Gamma _p}S(-d_{\gamma})v_{\gamma}, \qquad K_0=S.
 $$
\end{itemize} 
 Then an $S$-linear map $d:K_p \to K_{p-1}$ is defined by
$$
d(v_{\gamma})=\sum _{j=1}^p(-1)^{j+1}f_{i_j}v_{\gamma\setminus
	i_j}, \quad d(v_{(i_1)})=f_{i_1}.
$$
The Koszul complex $K_{\bul}$  is a graded complex with
  differentials of degree $0$, and the following is exact
(see 1.6.14 and 1.6.15 of \cite{c-m}):
$$
0 \to K_m  \xrightarrow{d} \cdots
 \xrightarrow{d}K_1   \xrightarrow{d} K_0 \to S/I \to 0,
$$
where $ K_0 \to S/I$ is the canonical surjection.
 Let $K_p^{\chi} \subset K_p$ be a submodule defined by
$$
K_p^{\chi}=\bigoplus_{\gamma \in
 \Gamma _p}S(-d_{\gamma})^{\chi
 \chi_{\gamma}^{-1}}v_{\gamma}; \quad K_0^{\chi}=S^{\chi}.
$$
After some computation we obtain
\begin{lem}
 \label{l:subcomplex}
We have $d(K_p^{\chi})\subset K_{p-1}^{\chi}$.
Moreover, the following complex is exact:
$$
0 \to K_m^{\chi}  \xrightarrow{d} \cdots
 \xrightarrow{d}K_1^{\chi} \xrightarrow{d} K_0^{\chi} \to (S/I)^{\chi} \to 0.
$$  
\end{lem}
If we  denote the  Hilbert series of a graded module $\cM$ by
$H_{\cM}(t)$, then \lemref{l:subcomplex} immediately implies 
$$
H^{\chi}(t)=\sum_{p=0}^m
	     (-1)^{p}H_{K_p^{\chi}}(t).
$$
Obviously,
$$
 H_{K_p^{\chi}}(t)=\sum _{\gamma
	     \in \Gamma
	     _p}H_{S(-d_{\gamma})^{\chi\chi_{\gamma}^{-1}}}(t), \quad
  H_{S(-d_{\gamma})^{\chi\chi_{\gamma}^{-1}}}(t)
  =t^ {d_{\gamma}}\cH^{{\chi\chi_{\gamma}^{-1}}}(t).
$$
By applying \thmref{t:molien}, we obtain the desired formula
$$
 H^{\chi}(t)=\frac{1}{|G|}\sum _{g \in G}
\frac{ \chi^{-1}(g)}{\prod_{j=1}^n(1-g_jt^{w_j})}
\sum_{p=0}^m	 
\sum _{\gamma	     \in \Gamma _p}
(-1)^{p}\chi_{\gamma}(g)t^{d_{\gamma}}.
$$

\subsection*{Proof of \proref{p:hp}}

We use the notation of \thmref{t:graded}.
For $w \in \cE$ (resp. $w \in \cN$), the $v$-degree of $z_w$
(resp. $\LF_v(f_{wj_w})$) is $m_{vw}$.
For $h \in \bH$, $h\cdot f_{wj_w}=\theta(h,E_w^*)f_{wj_w}$;
the same holds for $z_w$, $w \in \cE$
(see \S \ref{ss:H-action}).
Note that $g_j$  is determined by $g\cdot
z_j=g_jz_j$. 
Now, recall that $j_w$ moves from $1$ to $\delta_w-2$ if $w \in
\cN$ and $\delta_w-2=-1$ if $w \in \cE$.  Thus, it follows from
\thmref{t:molien2} that
$$
 \HP^{\chi}(t)=\frac{1}{|\bH|}\sum _{h\in  \bH}
\chi^{-1}(h) \prod_{w\in \cV}
\left(1-\theta(h, E^*_w)t^{m_{vw}}\right)^{\delta_w-2}.
$$


\providecommand{\bysame}{\leavevmode\hbox to3em{\hrulefill}\thinspace}
\providecommand{\MR}{\relax\ifhmode\unskip\space\fi MR }
\providecommand{\MRhref}[2]{%
  \href{http://www.ams.org/mathscinet-getitem?mr=#1}{#2}
}
\providecommand{\href}[2]{#2}

\end{document}